\documentclass[11pt]{smfart}

\usepackage{amsmath,a4wide}
\usepackage{amssymb}
\usepackage{xspace}
\usepackage{euscript}
\usepackage{mathrsfs} 
\usepackage[all]{xy}
\usepackage{hyperref}
\usepackage[francais]{babel}
\usepackage[T1]{fontenc}
\usepackage{textcomp}
\usepackage{smfthm}

\newtheorem{proposition}[subsection]{Proposition}	

\newtheorem{theoreme}[subsection]{Th\'eor\`eme}

\newtheorem{corollaire}[subsection]{Corollaire}

\theoremstyle{definition}

\newtheorem{remarque}[subsection]{Remarque}

\renewcommand{\theenumi}{\textbf{\roman{enumi}}} 

\renewcommand\theequation{\thesection.\textbf{\alph{equation}}} 

\makeatletter
\@addtoreset{equation}{subsection}
\makeatother

\renewcommand\thesubsubsection{\textbf{\thesubsection.\arabic{subsubsection}}}
\renewcommand\thesubsection{\textbf{\thesection.\arabic{subsection}}}

\newcommand{\Z}{\mathbb Z}

\newcommand{\R}{\mathbb R}
\newcommand{\C}{\mathbb C}
\newcommand{\F}{\mathbb F}

\newcommand{\Spec}{\mathrm{Spec}}

\renewcommand{\L}{\mathrm{L}}		

\def\commutatif{\ar@{}[rd]|{\circlearrowleft}}
\def\cartesien{\ar@{}[rd]|{\square}}

\renewcommand{\lim}{{\mathrm{lim}}} 

\def\ega#1#2{[{\bf \'EGA}~{\sc #1}~#2]} 				
\def\egazéro#1#2{[{\bf \'EGA}~$0_{\textsc{#1}}$~#2]}		

\numberwithin{equation}{subsubsection}
\renewcommand{\theequation}{\thesubsection.\Alph{equation}}

\newcommand{\isolong}{\buildrel{\sim}\over\longrightarrow}

\newcommand{\gsp}{\mathrm{GSp}_{2g}}

\newcommand{\agn}{\mathcal{A}_{g,n}}

\newcommand{\agnb}{\overline{\mathcal{A}}_{g,n}}

\title{Relèvement de formes modulaires de Siegel}
\author{Benoît Stroh}
\date{21 septembre 2009}
\email{benoit.stroh@gmail.com}
\address{Laboratoire Analyse, Géométrie et Applications, Université Paris 13, 93430 Villetaneuse, France}
\keywords{Siegel modular forms, toroïdal compactifications, Voronoï decomposition, Kodaira vanishing theorem, Deligne-Illusie Hodge theory}

\begin{document}
\maketitle




\begin{small}
{\textbf{Lifting Siegel modular forms. }}
In this note, we give explicit conditions under which cuspidal Siegel modular forms of genus~$2$ or~$3$ with coefficients in a finite field lift to cuspidal modular forms with coefficients in a ring of characteristic~$0$. This result extends a classical theorem proved by Katz for genus~$1$ modular forms. We use ampleness results due to Shepherd-Barron, Hulek and Sankaran, and vanishing theorems due to Deligne, Illusie, Raynaud, Esnault and Viehweg.
\end{small}

\bigskip

Dans cette note, nous donnons des conditions explicites sous lesquelles les formes modulaires de Siegel cuspidales de genre~$2$ ou $3$ à coefficients dans un corps fini se relèvent en des formes modulaires cuspidales à coefficients dans un anneau de caractéristique nulle. Nos résultats généralisent un théorème classique obtenu par Katz pour les formes de genre~$1$~\cite[th.~1.7.1]{Padic@Katz}. Nous utilisons des résultats d'amplitude de Shepherd-Barron~\cite{Perfect@Barron} et de Hulek et Sankaran~\cite{Nef@HulekSankaran}, ainsi que des théorèmes d'annulation dûs à Deligne, Illusie et Raynaud~\cite{Hodge@DeligneIllusie} et à Esnault et Viehweg~\cite{Annulation@EsnaultViehweg}. Je remercie le rapporteur pour ses corrections et pour m'avoir suggéré le corollaire~\ref{app1_coro3}.

\section{\'Enoncé des résultats}

Soient $g\geq 2$, $k\in\Z$ et $n\geq 3$ trois entiers et $M$ un $\Z[1/n]$-module. Notons $\agn$ l'espace de modules sur $\Spec(\Z[1/n])$ qui paramètre les schémas abéliens principalement polarisés de genre $g$ munis d'une structure de niveau principale en~$n$.
D'après~\cite[th.~IV.6.7]{Deg@FaltingsChai}, il existe une compactification toroïdale $\agnb$ de $\agn$ (elle dépend d'un choix combinatoire) et un schéma semi-abélien $G$ sur $\agnb$ qui étend le schéma abélien universel sur $\agn$.
Notons~$\omega$ le faisceau inversible sur $\agnb$ des formes volumes invariantes de $G$. Le principe de Köcher~\cite[prop.~V.1.8]{Deg@FaltingsChai} affirme que la restriction de $\agnb$ à $\agn$ induit un isomorphisme
$$\mathrm{H}^0\left(\agnb \, ,\,\omega^{k}\otimes M\right) \: \isolong \: \mathrm{H}^0\left(\agn \, ,\,\omega^{k}\otimes M\right) \: .$$
En particulier, le groupe $\mathrm{H}^0(\agnb \, ,\,\omega^{k}\otimes M)$ est indépendant du choix de $\agnb$~;~on appelle ses éléments les formes modulaires de genre $g$, de niveau $n$, de poids $k$ et à coefficients dans $M$.

Notons $D$ le complémentaire de $\agn$ dans $\agnb$ et munissons-le de sa structure schématique réduite. On dit qu'une forme modulaire est {cuspidale} si elle s'annule sur $D$. Notons
$$\mathrm{Cusp}(n,k,M) \: = \:  \mathrm{H}^0\left(\agnb  \, ,\, \omega^{k}(-D) \otimes M\right)$$ 
le groupe des formes cuspidales. D'après~\cite[prop.~V.1.9]{Deg@FaltingsChai}, il ne dépend pas non plus du choix de~$\agnb$. Si $M$ est un anneau, on a
$$\mathrm{Cusp}(n,k,M) \: = \:  \mathrm{H}^0\left(\agnb \times \Spec(M) \, ,\, \omega^{k}(-D)\right)\: .$$
\'Enonçons à présent le théorème principal de cette note.
\begin{theoreme} \label{app1_th1} Supposons $n\geq 3$, $g=2$ ou~$3$ et $k>g+1$. Soit~$M$ un~$\Z[\frac{1}{n}]$-module. Si~$g=2$, supposons~$M$ sans $2$-torsion, et si~$g=3$ supposons~$M$ sans $30$-torsion. Le morphisme naturel de changement de base induit un isomorphisme
$$\mathrm{Cusp}(n,k,\Z[1/n])\:\otimes_{\Z[1/n]} \: M \: \isolong \: \mathrm{Cusp}(n,k,M) \: .$$
\end{theoreme}

Ce théorème sera démontré dans la troisième partie. Notons qu'on peut s'affranchir de l'hypothèse~$n\geq 3$ en introduisant le niveau auxiliaire~$3n$ et en considérant les formes modulaires de niveau~$3n$ qui sont invariantes par le groupe~$\gsp(\Z/3\Z)$.

\begin{corollaire}\label{app1_coro2} Soient $n\geq 3$, $g=2$ ou $3$, et $k>g+1$ trois entiers. Soit~$\mathcal{O}$ l'anneau d'entiers d'un corps de nombres et $\mathfrak{m}$ un idéal maximal de~$\mathcal{O}$ ; notons $\kappa=\mathcal{O}/\mathfrak{m}$ et supposons que la caractéristique~$p$ de~$\kappa$ ne divise pas~$n$. Si $g=2$, supposons $p>2$ et si $g=3$, supposons $p>5$.
Le morphisme de changement de base induit un isomorphisme
$$\mathrm{Cusp}(n,k,\mathcal{O}[1/n]) \: \otimes_{\mathcal{O}[1/n]} \: \kappa  \:  \isolong \:  \mathrm{Cusp}(n,k,\kappa)\: .$$
\end{corollaire}

Il suffit en effet d'appliquer successivement le théorème~\ref{app1_th1} à~$M=\mathcal{O}[1/n]$ et~$M=\kappa$ pour démonter le corollaire~\ref{app1_coro2}.

\begin{corollaire} \label{app1_coro3} Soient $n\geq 3$, $g=2$ ou $3$, et $k>g+1$ trois entiers. Soit~$\mathcal{O}$ un anneau de valuation discrète d'idéal maximal~$\mathfrak{m}$ ; notons $\kappa=\mathcal{O}/\mathfrak{m}$ et supposons que la caractéristique~$p$ de~$\kappa$ ne divise pas~$n$. Si $g=2$, supposons $p>2$ et si $g=3$, supposons $p>5$. 
Notons $\mathcal{H}$ la sous-algèbre commutative de $\mathrm{End}(\mathrm{Cusp}(n,k,\mathcal{O}))$ engendrée par les opérateurs de Hecke de niveau premier à~$n$. Soit~$\bar{f}\in \mathrm{Cusp}(n,k,\kappa)$ propre pour~$\mathcal{H}$ de valeur propre généralisée~$\bar{\chi}:\mathcal{H} \rightarrow \kappa$. Il existe~$f' \in \mathrm{Cusp}(n,k,\mathcal{O})$ propre pour~$\mathcal{H}$ de valeur propre généralisée~${\chi}:\mathcal{H} \rightarrow \mathcal{O}$ telle que~$\chi\: \mathrm{mod}\: \mathfrak{m} = \bar{\chi}$.
\end{corollaire}

Le corollaire~\ref{app1_coro3} affirme que toute valeur propre généralisée de~$\mathcal{H}$ dans~$\mathrm{Cusp}(n,k,\kappa)$ se relève en une valeur propre généralisée dans~$\mathrm{Cusp}(n,k,\mathcal{O})$. Par contre, il ne dit rien sur le relèvement des vecteurs propres, \textit{ie}. ne garantit pas que~$f' \: \mathrm{mod} \: \mathfrak{m} = \bar{f}$.
Pour démontrer le corollaire~\ref{app1_coro3}, on procède en deux temps : on utilise le théorème~\ref{app1_th1} pour prouver que~$\mathrm{Cusp}(n,k,\kappa)=\mathrm{Cusp}(n,k,\mathcal{O})\otimes\kappa$, puis l'on applique le lemme de Deligne-Serre~\cite[lem.~6.11]{Mod@DeligneSerre}.
 
\begin{remarque} Le théorème~\ref{app1_th1} et les corollaires~\ref{app1_coro2} et~\ref{app1_coro3} concernent uniquement les formes cuspidales.  Nous ne pensons pas que nos méthodes puissent s'adapter au cas des formes non nécessairement cuspidales.
\end{remarque}

Nous prouverons également la proposition suivante, qui affirme qu'il n'existe pas de forme modulaire (cuspidale ou non) non nulle de genre $2$ ou~$3$, de poids strictement négatif et à coefficients dans un corps de caractéristique~$>5$ (\textit{cf}. prop. \ref{prop_negatif}).

\begin{proposition} Si $g=2$ ou $3$, $n\geq 3$ et $k<0$, on a $\mathrm{H}^0(\mathcal{A}_{g,\, n}\times \kappa\, ,\omega^{k})=0$ pour tout corps~$\kappa$ de caractéristique~$>5$.
\end{proposition}

\begin{remarque} Soit~$p$ un nombre premier. Le corollaire~\ref{app1_coro2} permet de préciser quantitavement un théorème de classicité des formes modulaires de Siegel $p$-adiques  ordinaires~\cite[th.~1.1]{Padic@Hida}. Soit~$f$ une forme modulaire de Siegel $p$-adique~\cite[déf.~1.4.3]{These@Pilloni} ordinaire, cuspidale, de poids classique et parallèle~$k\in\mathbb{Z}$, de niveau~$n\geq 3$, et de genre~$g$ égal à $2$ ou~$3$. Rappelons que~$f$ est appelée \og classique \fg si elle provient d'un élément de~$\mathrm{Cusp}(n,k,\Z_p)$. D'après~\cite[th.~1.1.(4)]{Padic@Hida} ou~\cite[th.~1.6.6]{These@Pilloni}, $f$ est classique lorsque~$k$ est supérieur à une constante indéterminée dépendant de~$p$, de~$g$ et de~$n$. Le corollaire~\ref{app1_coro2} permet de prouver que l'estimation $k>g+1$ suffit à assurer la classicité de~$f$ si~$p>5$. Il faut en effet reprendre la démonstration de~\cite[th.~1.6.6]{These@Pilloni} et se rendre compte que l'hypothèse \og $k$ assez grand \fg n'est utilisée que dans~\cite[th.~1.6.1]{These@Pilloni}. Mais dans le cas présent,~\cite[th.~1.6.1]{These@Pilloni} est équivalent au corollaire~\ref{app1_coro2}, qui est valable si~$k>g+1$ et~$p>5$.
\end{remarque}

\section{Amplitude et annulation} 

Dans cette partie, nous rappelons brièvement les résultats utilisés dans la démonstration du théorème~\ref{app1_th1}. On suppose~$g\geq 2$ et~$n\geq 3$ dans toute cette partie.

\subsection{Compactifications toroïdales}
Notons $B$ le $\Z$-module des formes bilinéaires symétriques entières sur $\Z^g$ et $C$ le cône de $B\otimes \R$ des formes bilinéaires semi-définies positives à radical rationnel. Les compactifications toroïdales de $\agn$ sont associées à des décompositions polyédrales admissibles $\mathrm{GL}_g(\Z)$-équivariantes de $C$~\cite[déf.~IV.2.2]{Deg@FaltingsChai}. Considérons une telle décomposition $\Sigma$ et notons $\agnb$ la compactification toroïdale qui lui est associée~\cite[th. IV.6.7]{Deg@FaltingsChai}. L'espace algébrique $\agnb$ est propre sur $\Spec(\Z[1/n])$, lisse sur $\Spec(\Z[1/n])$ si et seulement $\Sigma$ est lisse~\cite[déf.~IV.2.3]{Deg@FaltingsChai}, et est un schéma projectif si $\Sigma$ est polarisée~\cite[déf.~V.2.4 et~th.~V.5.8]{Deg@FaltingsChai}.
Notons~$D$ le complémentaire de $\agn$ dans $\agnb$ et munissons-le de sa structure de schéma réduit. Il définit un diviseur de Cartier de $\agnb$ qui est à croisements normaux si $\Sigma$ est lisse. Notons~$\{f_i\}$ un système d'équations locales de $D$ dans $\agnb$. Comme les éventuelles singularités de $\agnb$ sont toriques, le faisceau
$$\Omega^1_{\agnb}\left(\mathrm{log}\: D\right)\: = \: \Omega^1_{\agnb} \left( \frac{\mathrm{d} f_i}{f_i} \right)$$
est localement libre. D'après~\cite[th. IV.6.7]{Deg@FaltingsChai}, le schéma abélien universel sur~$\agn$ s'étend en un schéma semi-abélien~$G$ sur~$\agnb$. Notons~$e$ la section neutre de~$G$. 
Toujours d'après~\cite{Deg@FaltingsChai}, l'application de Kodaira-Spencer induit un isomorphisme
$$\mathrm{Sym}^2(e^*  \Omega^1_{G/\agnb} )\: \isolong \: \Omega^1_{\agnb}\left(\mathrm{log}\: D\right).$$
Notons $\mathcal{K}=\mathrm{det}\left(\Omega^1_{\agnb}\right)$ le faisceau dualisant de $\agnb$. Puisque~$\omega=\mathrm{det}(e^*  \Omega^1_{G/\agnb})$, l'isomorphisme de Kodaira-Spencer induit un isomorphisme $\omega^{g+1} \isolong \mathcal{K}(D)$~\cite[p.225 (d)]{Deg@FaltingsChai}.

\subsection{Décompositions de Voronoï} 
Voronoï a construit deux décompositions polyédrales admissibles de $C$~\cite{Parfait@Voronoi}. L'une d'entre elles est appelée \textit{décomposition du cône parfait} ou \textit{première décomposition de Voronoï}, et l'autre \textit{seconde décomposition de Voronoï}. Igusa a construit une troisième décomposition de $C$, dite \textit{du cône central}~\cite{Nef@HulekSankaran}. Ces trois décompositions coïncident si $g=2$ ou $3$~\cite[p.~660]{Nef@HulekSankaran}. De plus, elles sont lisses si $g=2$ ou $3$~\cite[p.~255 et~256]{Nef@Hulek}. Si $g=4$, la seconde décomposition de Voronoï est lisse~\cite[p.~661]{Nef@HulekSankaran} et raffine la décomposition du cône parfait~\cite[p.~660]{Nef@HulekSankaran}.

\subsection{Amplitude du fibré de Hodge}

Notons $\overline{\mathcal{A}}_{g,\,n}^{\mathrm{{parf}}}$ la compactification toroïdale de~$\agn$ associée à la décomposition du cône parfait et $D^{\mathrm{parf}}$ son bord réduit. Le schéma $$\overline{\mathcal{A}}_{g,\,n}^{\mathrm{{parf}}}$$ est lisse si $g=2$ ou~$3$, et $D^{\mathrm{parf}}$ est un diviseur à croisements normaux. Soit $a$ un entier relatif. D'après~\cite[th.~4.1]{Perfect@Barron}, le faisceau $\omega^{a}(-D^{\mathrm{{parf}}})$ est ample sur $$\overline{\mathcal{A}}_{g,\,n}^{\mathrm{{parf}}}$$ si et seulement si $a>12/n$.
Remarquons que ce résultat d'amplitude a d'abord été établi par Hulek et Sankaran~\cite{Nef@HulekSankaran} dans le cas particulier $g=2$ ou $3$ considéré dans la suite de ce texte. Cependant, contrairement à Shepherd-Barron, ils n'énoncent leur théorème que sur~$\Spec(\C)$.

\subsection{Théorème de Kodaira}

Soit $\kappa$ un corps parfait de caractéristique $p>0$. Donnons-nous un schéma $X$ propre et lisse purement de dimension $d\leq p$ sur $\Spec(\kappa)$, un diviseur à croisements normaux $D$ de $X$, et un faisceau inversible $\mathcal{L}$ sur $X$. Notons $W_2(\kappa)$ l'anneau des vecteurs de Witt de longueur~$2$ de $\kappa$ et supposons que $X$ et $D$ se relèvent en~$\tilde{X}$ et $\tilde{D}$ sur $\Spec(W_2(\kappa))$, où~$\tilde{X}$ est lisse sur~$\Spec(W_2(\kappa))$ et $\tilde{D}$ est un diviseur à croisements normaux de~$\tilde{X}$.


Supposons également qu'il existe un entier $\nu_0$ tel que $\mathcal{L}^{\nu}(-D)$ soit ample pour tout $\nu>\nu_0$. D'après~\cite[prop.~11.5]{Annulation@EsnaultViehweg}, on~a 
\begin{equation} \label{app1_eq1}
\mathrm{H}^j\left(X,\:\Omega_{X/\kappa}^i(\mathrm{log}\: D )  \otimes \mathcal{L}^{-1} \right) \: = \: 0
\end{equation}
pour $i+j<d$. En posant~$i=0$ dans l'égalité~(\ref{app1_eq1}) et en utilisant le fait que~$\Omega_{X/\kappa}^0(\mathrm{log}\: D )=\mathcal{O}_X$, on trouve en particulier
$$\mathrm{H}^j\left(X,\:\mathcal{L}^{-1} \right) \: = \: 0 $$
pour tout~$j<d$.
On en déduit par dualité de Serre que
\begin{equation} \label{app1_eq2}
\mathrm{H}^j(X,\:\mathcal{K}_{X/\kappa}\otimes\mathcal{L})\:=\:0
\end{equation}
pour tout~$j>0$, où $\mathcal{K}_{X/\kappa}=\Omega^d_{X/\kappa}$ désigne le faisceau dualisant de $X$. 

\begin{remarque} L'assertion~(\ref{app1_eq1}) est une version faible mais valable sur $\kappa$ du théorème de Kawamata-Viehweg pour les diviseurs nefs~\cite[coro.~5.12.c]{Annulation@EsnaultViehweg}, qui n'est démontré pour l'instant que sur $\C$. Le théorème de Kawamata-Viehweg sur $\kappa$ résulterait de l'assertion~(\ref{app1_eq1}) et d'une éventuelle résolution des singularités des schémas de type fini sur $W_2(\kappa)$. Nous renvoyons à~\cite[rem.~11.6.b]{Annulation@EsnaultViehweg} pour plus de détails.
\end{remarque}

\section{Démonstration du théorème~\ref{app1_th1}}

Nous nous plaçons à présent sous les hypothèses du théorème~\ref{app1_th1}. En particulier, $g=2$ ou~$3$, $n\geq 3$ et~$k>g+1$. Comme la formation de la cohomologie cohérente d'un schéma quasi-cohérent commute aux limites inductives, il suffit de démontrer le théorème~\ref{app1_th1} pour $M$ de type fini sur $\Z[1/n]$, puis pour $M=\F_p$ avec $p>2$ si $g=2$ et $p>5$ si $g=3$. Ainsi, il suffit de prouver que le morphisme de changement de base induit un isomorphisme
$$ \mathrm{H}^0\left(\:\overline{\mathcal{A}}_{g,n}^{\mathrm{{parf}}}\, ,\, \omega^{k}(-D^{\mathrm{parf}})\right)\otimes \F_p \: \isolong \: \mathrm{H}^0\left(\:\overline{\mathcal{A}}_{g,n}^{\mathrm{{parf}}}\times \Spec(\F_p)\, ,\, \omega^{k}(-D^{\mathrm{parf}})\right)\: .$$
D'après~\ega{3}{7.5.3 et 7.7}, il suffit de démontrer que 
$$\mathrm{H}^1\left(\:\overline{\mathcal{A}}_{g,\,n}^{\mathrm{{parf}}}\times \Spec(\F_p)\, ,\, \omega^{k}(-D^{\mathrm{parf}})\right) \: = \: 0\: .$$
Le schéma  $\overline{\mathcal{A}}_{g,\,n}^{\mathrm{{parf}}}\times \Spec(\F_p)$
est propre et lisse sur $\Spec(\F_p)$ car $g=2$ ou $3$. Il est purement de dimension $g(g+1)/2$ et se relève en un schéma lisse sur~$\Spec(\Z/p^2\Z)$. De même, $D^{\mathrm{parf}}$ est un diviseur à croisements normaux qui se relève à $\Spec(\Z/p^2\Z)$. Notons $\mathcal{L}=\omega^{k-1-g}$ ; c'est un faisceau inversible sur $$\overline{\mathcal{A}}_{g,\,n}^{\mathrm{{parf}}}\times \Spec(\F_p)\: .$$
D'après~\cite[th.~4.1]{Perfect@Barron}, $\mathcal{L}^{\nu}(-D^{\mathrm{parf}})$ est ample sur  $\overline{\mathcal{A}}_{g,\,n}^{\mathrm{{parf}}}$ pour tout $\nu > 12/n(k-1-g)$.
D'après~(\ref{app1_eq2}), on a
$$\mathrm{H}^1\left(\:\overline{\mathcal{A}}_{g,\,n}^{\mathrm{{parf}}}\times \Spec(\F_p)\:, \:\mathcal{K}\otimes \mathcal{L} \right) \:=\: 0$$
car $p\geq g(g+1)/2$ par hypothèse. On conclut la démonstration du théorème~\ref{app1_th1} en remarquant que $$\mathcal{K}\otimes \mathcal{L} = \omega^{k}(-D^{\mathrm{parf}})\: .$$

Montrons à présent qu'il n'existe pas de forme modulaire de Siegel de genre~$2$ ou $3$, de poids~$<0$ et à coefficients dans un corps~$\kappa$ (non nécessairement parfait) de caractéristique $p>5$. Il nous faut donc voir que~$\mathrm{H}^0(\mathcal{A}_{g,\, n}\, ,\omega^{k}\otimes \kappa)=0$. On raisonne comme au début de la partie~3 et on se ramène à prouver que~$\mathrm{H}^0(\mathcal{A}_{g,\, n}\times \F_p\, ,\omega^{k})=0$. Il suffit de poser $$X=\overline{\mathcal{A}}_{g,\,n}^{\mathrm{parf}}\times \F_p\: ,$$ $i=0$, $j=0$ et $\mathcal{L}=\omega^{-k}$ dans~(\ref{app1_eq1}) pour conclure, c'est-à-dire démontrer le résultat suivant.

\begin{proposition} \label{prop_negatif} Si $g=2$ ou $3$, $n\geq 3$ et $k<0$, on a $\mathrm{H}^0(\mathcal{A}_{g,\, n}\times \kappa\, ,\omega^{k})=0$ pour tout corps~$\kappa$ de caractéristique~$>5$.
\end{proposition}

\section{Cas du genre quatre}

On aimerait prouver un énoncé de changement de base analogue au théorème~\ref{app1_th1} dans le cas où $g=4$. Comme la décomposition du cône parfait n'est plus lisse, on ne peut pas appliquer le théorème d'annulation de Kodaira à $\overline{\mathcal{A}}_{4,\,n}^{\mathrm{{parf}}}$. En revanche, la seconde décomposition de Voronoï est lisse et raffine la décomposition du cône parfait d'après~\cite[p.~661]{Nef@HulekSankaran}. Notons $$\overline{\mathcal{A}}_{4,\,n}^{\mathrm{{Vor}}}$$ la compactification toroïdale associée à la seconde décomposition de Voronoï. Il existe un morphisme propre et birationnel
$$\pi \, : \, \overline{\mathcal{A}}_{4,\,n}^{\mathrm{{Vor}}} \: \longrightarrow \: \overline{\mathcal{A}}_{4,\,n}^{\mathrm{{parf}}} $$
induit par les raffinements entre décompositions.
Notons $E$ le diviseur exceptionnel de $\pi$, qui est de support inclus dans~$D^{\mathrm{{Vor}}}$. En combinant le théorème~1.8, la remarque~1.9 et le corollaire~1.11 de~\cite{Nef@HulekSankaran}, on peut montrer que 
$\omega^a(-b D^{Vor} - c E)$ est ample sur $$\overline{\mathcal{A}}_{4,\,n}^{\mathrm{{Vor}}}$$ si et seulement si $a>\frac{12b}{n}$ et $c>4b>\frac{8}{9}c$. Par exemple, si $a>36/n$ alors $$\omega^a(-3D^{\mathrm{Vor}} - 13 E)$$ est ample sur $\overline{\mathcal{A}}_{4,\,n}^{\mathrm{{Vor}}}$. On peut appliquer le théorème de Kodaira à la variété propre et lisse $$\overline{\mathcal{A}}_{4,\,n}^{\mathrm{{Vor}}}\times\Spec(\F_p)$$ pour relever en caractéristique nulle des formes de Siegel de genre $4$, de niveau~$>36$, de poids~$>5$, à coefficients dans $\F_p$, qui s'annulent suffisamment sur $D^{\mathrm{Vor}}$. L'ordre d'annulation requis est au moins~$2$ et on peut le déterminer en calculant explicitement~$E$.


\providecommand{\bysame}{\leavevmode ---\ }
\providecommand{\og}{}
\providecommand{\fg}{}
\providecommand{\smfandname}{et}
\providecommand{\smfedsname}{\'eds.}
\providecommand{\smfedname}{\'ed.}
\providecommand{\smfmastersthesisname}{M\'emoire}
\providecommand{\smfphdthesisname}{Th\`ese}

\end{document}